\magnification\magstep1
\baselineskip=18pt
%lines for a solid square at end of proofs
\def \qed {\vrule height6pt  width6pt depth0pt}

\def \Bbb {\bf}
\def \dim {{\rm dim\,}}
\def \det {{\rm det\,}}
\def \Log{{\rm Log \,}}
\def \rk{{\rm rk \,}}
\def \tr{{\rm tr \,}}
\def \Re{{\rm Re \,}}
\centerline{{\bf The proportional UAP characterizes weak 
Hilbert
spaces}\footnote*{Both authors were supported in part by NSF 
DMS 87-03815}}
\vskip12pt
\centerline {by W. B. Johnson and G. Pisier}\vskip40pt

{\bf Abstract:}  We prove that a Banach space has the uniform 
approximation
property with proportional growth of the uniformity function 
iff it is a
weak Hilbert space.

\vfill\eject

\noindent {\bf Introduction}

\noindent The ``weak Hilbert spaces'' were introduced and 
studied in [P 2]. 
Among the many equivalent characterizations in [P 2] perhaps 
the simplest
definition is the following.  A Banach space is a weak Hilbert 
space if there is
a constant C such that for all $n$, for all $n$-tuples 
$(x_1,\cdots , x_n)$ and
$(x^*_1,\cdots ,x^*_n)$ in the unit balls of $X$ and $X^*$ 
respectively, we have

$${|{\rm det}(<x^*_i,x_j>)|}^{1\over n} \le C.$$

\noindent The first example of a non Hilbertian weak Hilbert 
space was obtained
by the first author (cf. [FLM], Example 5.3 and [J]).\vskip6pt

\noindent Recall that a Banach space X has the uniform 
approximation property
(in short UAP) if there is a constant $K$ and a function 
$n\rightarrow f(n)$
such that for all $n$ and all $n$-dimensional subspaces 
$E\subset X$, there is
an operator $T: X\rightarrow X$ with $\rk(T)\le f(n)$ such that 
$\|T\|\le K$ and
$T_{|E} = I_{|E}$.\vskip6pt

\noindent For later use, given $K> 1$ we introduce

$$k_X(K,n) = \sup_{\scriptstyle E\subset X\atop\scriptstyle
\dim E=n }\ \ \inf\ \{\rk(T)\}$$

\noindent where the infimum runs over all $T: X\rightarrow X$ 
such that $\|T\|
\le  K$ and $T_{|E} = I_{|E}$.  

\noindent Note that $X$ has the UAP iff there is a constant $K$ 
such that
$k_X(K,n)$ is finite for all $n$; we then say that $X$ has the 
$K$-UAP.  The asymptotic growth of the function
$n\rightarrow k_X(K,n)$ provides a quantitative measure of 
the  UAP of the space $X$.

\noindent For instance, if $X$ is a Hilbert space we have clearly 
$k(1,n) = n$,
hence if $X$ is isomorphic to a Hilbert space there is a constant 
$K$ such that

$$k_X(K,n) = n\quad {\rm for\ all}\ n.$$

\noindent The converse is also true by the complemented 
subspace theorem of
Lindenstrauss-\break Tzafriri [LT 1].\vskip6pt

\noindent The main result in this paper can be viewed as an 
analogous statement
for weak Hilbert spaces, as follows.\vskip12pt

\proclaim Main Theorem. {\sl A Banach space $X$ is a weak 
Hilbert space
iff there are constants $K$ and $C$ such that
$$k_X(K,n)\le Cn\quad {\rm for\ all}\ n.\leqno (0.1)$$}

\noindent That is, proportional asymptotic behavior of the 
uniformity function
in the definition of the UAP characterizes weak Hilbert 
spaces.\vskip6pt

\noindent It was proved in [P 2] that weak Hilbert spaces 
have 
the UAP
but no estimate of the function $n\rightarrow k_X(K,n)$ was 
obtained.\vskip6pt

\noindent For the purposes of this paper we will say that $X$ 
has the
proportional UAP if there are constants $K$ and $C$ such that 
(0.1) holds.\vskip6pt
\noindent The authors thank V. Mascioni and G. 
Schechtman for 
several discussions concerning the material in this paper.
\vskip12pt

\noindent {\bf $\S 1.$ Weak Hilbert spaces have proportional 
UAP}\vskip6pt

\noindent We first recall a characterization of weak Hilbert 
spaces in terms of
nuclear operators.  Recall that an operator $u:X\rightarrow X$ 
is called
nuclear if it can be written

$$u(x) = \sum^\infty_{n=1} x^*_n(x)x_n$$

\noindent with $x^*_n \in X^*, x_n\in X$ such that $\sum 
\|x^*_n\|\  \|x_n\| <
\infty$.  Moreover the nuclear norm $N(u)$ is defined as

$$N(u) = \inf\ \Big\{\sum \|x^*_n\|\  \|x_n\| \Big\}$$

\noindent where the infimum runs over all possible 
representations.  We also
recall the notation for the approximation numbers

$$\forall k\ge 1\quad a_k(u) = \inf\ \{\|u-v\|\ \ |\  v: 
X\rightarrow X,\ \ \rk(v)<k\}.$$

\noindent By [P 2], a Banach space $X$ is a weak Hilbert space 
iff 
there is a
constant $C$ such that for all nuclear operators $u: 
X\rightarrow X$ we have

$$\sup_{k\ge 1}\ k a_k (u) \le C N(u)\leqno(1.1)$$

\noindent The following observation is identical to reasoning 
already used by
V. Mascioni [Ma 2].

\proclaim  Proposition 1.1.  {\sl Let $X$ be a weak Hilbert 
space.  Assume that
there is a constant $K'$ such that for all $n$ and all $n$ 
dimensional
subspaces $E\subset X$ there is an operator $u: X\rightarrow 
X$ such that
$u_{|E} = I_{|E}$, $\|u\| \le K'$ and $N(u) \le K'n$.  Then $X$ has 
the
proportional UAP.  (Recall that if $u$ has finite rank then $N(u) 
\le
\rk(u)\|u\|$, hence the converse to the preceding implication is 
obvious.)}

\noindent {\bf Proof}:  Let $u$ be as in the preceding 
statement.  We use (1.1)
with $k = [2CK'n] + 1,$ so that
$$a_k(u) \le CN(u) k^{-1} \le CK'nk^{-1}\le {1\over 2}.$$

\noindent This means that there is an operator $v: 
X\rightarrow X$ with
$\rk(v)\le 2CK'n$ such that $\|u-v\|\le {1\over 2}$.  By 
perturbation, it follows that the operator
$$V = v-u+I$$

\noindent is invertible on $X$ with $\|V^{-1}\|\le 2$.  Moreover 
we have
$$V_{|E} = v_{|E}.\leqno (1.2)$$

\noindent It follows that if we let $T = V^{-1}v$, then we have
$$\|T\| \le \|V^{-1}\|\ \|v\| \le 2 (\|u\| + \|u-v\|) \le 2K' + 1,$$

\noindent also $\rk(T) \le \rk(v)\le 2CK'n$ and  $T_{|E} = 
I_{|E}$ 
by (1.2).

\noindent We conclude that $X$ has the UAP with $k_X(K,n)\le
2CK'n$, where $K=2K' + 1$. \quad\qed

\noindent We will use duality via the following proposition (a 
similar kind of criterion was used by Szankowski [S] to prove 
that certain spaces {\sl fail} the UAP):\vskip6pt

\proclaim Proposition 1.2.  {\sl Let $X$ be a reflexive Banach 
space with the approximation property; in short, AP; 
let $\alpha ,\beta$ be positive constants; and let $n\ge 1$ 
be an
integer.  The following are equivalent.
\item {(i)}  For all nuclear operators $T_1$, $T_2$ on $X$ such 
that $T_1 + T_2$
has rank $\le n$, we have
$$|\tr(T_1 + T_2)|\le \alpha N(T_1) + \beta n \|T_2\|.$$
\item {(ii)}  Same as (i) with $T_1$, $T_2$ of finite rank.
\item {(iii)}  For any subspace $E\subset X$ with dimension 
$\le n$ there is an
operator $u: X\rightarrow X$ such that $u_{|E} = I_E,\ \|u\| \le 
\alpha$ and
$N(u) \le \beta n$.}

\noindent {\bf Proof}:  (i) $\Rightarrow$ (ii) is trivial.

\noindent Assume that (ii) holds.

\noindent We equip $X^*\otimes X$ with the norm $|w| = \inf\ 
\{\alpha N(T_1) +
\beta n\|T_2\|\}$ where the infimum runs over all 
decompositions
$$u = T_1 + T_2 \ {\rm with}\ T_1\ {\rm and}\ T_2\ {\rm in}\ 
X^* \otimes X$$

\noindent (identified with the set of finite rank operators on 
$X$).  On
$X^*\otimes X$ this norm is clearly equivalent to the operator 
norm.

\noindent Now let $E\subset X$ be a fixed subspace with $\dim 
(E) \le n$.  Let
${\cal S} \subset X^*\otimes X$ be the subspace $X^*\otimes E$ 
of all the
operators on $X$ with range in $E$.  On this linear subspace the 
linear form $\xi$
defined by $\xi(w) = \tr(w)$ has norm $\le 1$ relative to 
$|\cdot 
|$ by our
assumption (ii).\vskip6pt

\noindent Therefore there is a Hahn-Banach extension 
${\tilde\xi}$ defined on
the whole of $X^*\otimes X$ which extends $\xi$ and satisfies

$$|{\tilde\xi}(w)| \le |w|\quad \forall w\in X^*\otimes X.\leqno 
(1.3)$$

\noindent By classical results, ${\tilde\xi}$ can be identified 
with an
integral operator $u: X\rightarrow X^{**}$.  Since $X$ is 
reflexive, $u$ is
actually a nuclear operator on $X$, and we have ${\tilde\xi}(w) 
= \tr(wu)$ for
all $w$ in $X^*\otimes X$.

\noindent Since ${\tilde\xi}$ extends $\xi$, we must have

$$\forall x^*\in X^*\quad \forall e\in E\quad <{\tilde\xi},\ 
x^*\otimes
e> = \tr(x^*\otimes e)=x^*(e)\quad {\rm hence}\ x^*(ue) = 
x^*(e).$$

\noindent Equivalently
$$u_{|E} = I_{|E}.$$

\noindent On the other hand, by (1.3) we have
$$|\tr(u T_1)| \le \alpha N(T_1)\quad {\rm and}\ |\tr (u T_2)\ 
|\ 
\le
\beta n \| T_2\|$$

\noindent for all finite rank operators $T_1$ and $T_2$ on $X$.

\noindent This implies $\|u\|\le \alpha$ and (again using the 
reflexivity of
$X$) $N(u) \le \beta n$.

\noindent This shows that (ii)$\Rightarrow$(iii).

\noindent Finally we show (iii)$\Rightarrow$(i).  Assume (iii).  
Let $T_1, T_2$
be as in (i), let $E$ be the range of $T_1 + T_2$ and let $u$ be 
as in (iii).
Then we have
$T_1 + T_2 = u(T_1 + T_2)$ 
hence since $X$ has the AP (which ensures
that $|\tr(T)|\le N(T)$ for all nuclear operator $T: X\rightarrow 
X$) we have
$|\tr(T_1+T_2)| = |\tr (uT_1) + \tr (uT_2)|$

$$\eqalign {&\le \|u\| N(T_1) + N(u) \|T_2\|\cr&\le\alpha (T_1) 
+ \beta n
\|T_2\|.\quad \qed }$$

\noindent {\bf Remark}:  Note that (i) is also equivalent to (i'):
\vskip6pt  
(i') For all
$T_1, T_2$ on $X$ such that $(T_1 + T_2)$ has rank $\le n$, we 
have
$$N(T_1+T_2)\le \alpha N(T_1) + \beta n \|T_2\|.$$
\noindent Indeed; (assuming the AP and reflexivity) we have
$$N(T_1 + T_2) = \sup \{\tr(S(T_1 + T_2))\ ;\ S: X\rightarrow X, 
\ \|S\|\le 1\}.$$

\noindent This shows that (i)$\Rightarrow$ (i').  Since $X$ has 
the AP
the converse is obvious.

\noindent Of course, a similar remark holds for (ii).

\noindent {\bf Remark}:  If $X$ is not assumed to have the AP 
a variant of
Proposition 1.2 will still hold provided we use the projective 
tensor norm on
$X^*\otimes X$ instead of the nuclear norm.\vskip6pt

\noindent We will use the following result already exploited in 
[P 2] to prove
that weak Hilbert spaces have the AP.  Whenever $u: 
X\rightarrow X$ is
a finite rank operator, we denote by ${\rm det}\,(I+u)$ the 
quantity

$$\Pi (1 + \lambda_j(u))$$

\noindent where $\{\lambda_j(u)\}$ are the eigenvalues of $u$ 
repeated
according to their algebraic multiplicity.  Equivalently, ${\rm 
det}\,(I+u)$ is equal
to the ordinary determinant of the operator $(I + u)_{|E}$ 
restricted to any
finite dimensional subspace $E\subset X$ containing the range 
of $u$.\vskip6pt

\proclaim Lemma 1.3.  {\sl Let $u, v$ be finite rank operators 
on a weak
Hilbert space $X$ with $\rk(u)\le n$.
Then we have $$|{\rm det\ }(I+u+v)| \le \Big (\sum^n_{j=0} 
{C^j\over j!} N(u)^j\Big )
\exp CN(v)\leqno (1.4)$$
where $C$ is the ``weak Hilbert space constant'' of $X$; that is 
to
say,
$$C=\sup_{\scriptstyle x_i\in B_X\atop \scriptstyle x^*_i \in 
B^*_X}
{| \det (<x^*_i, x_j>)|}^{1/n}.\leqno (1.5)$$}

\noindent For the proof we refer to [P 3] p. 229.  Note that if 
$N(u) \ge 1$
then (1.4) implies for all complex numbers $z$,
$${\rm det\ }(I + z(u+v))\le N(u)^n \exp\{C|z| + C|z| 
N(v)\}.\leqno 
(1.6)$$
\noindent Let $f(z) = \det (I + z(u+v))$. Then $f$ is a 
polynomial function of
$z\in {\Bbb C}$ such that 

$$f(0) = 1 \ {\rm and}\ f'(0) = \tr(u+v).$$

\noindent In [G], Grothendieck showed that the function 
$u\rightarrow {\rm det}\,(I+u)$
is uniformly continuous on $X^*\otimes X$ equipped with the 
projective norm,
and therefore extends to the completion $X^*{\hat\otimes}X$.  
This shows that
if $X$ has the AP, the determinant ${\rm det}\,(I+v)$ can be 
defined unambiguously for
any nuclear operator $v$ on $X$.  As shown in [G], the function 
$z\rightarrow
{\rm det}\,(I+z(u+v))$ is an entire function satisfying (1.4) if 
$u$ is 
of rank $\le n$ and
$v$ possibly of infinite rank.  We use this extension in Theorem 
1.5 below, but
in the proof of our main result the special case of $v$ of finite 
rank in
Theorem 1.5 is sufficient.  This makes our proof more 
elementary.

\noindent We will make crucial use of the following classical 
inequality of
Carath\'eodory; we include the proof for the convenience of the 
reader.\vskip6pt

\proclaim Lemma 1.4.  {\sl Let $h$ be an analytic function in a 
disc $D_R =
\{z\in {\Bbb C}\ ;\quad |z|<R\}$ such that $h(0)=0$.  Then for 
any $0<r<R$ we
have 
$$|h'(0)|\le {2\over r} \sup_{|z|=r} \Re(h(z)).$$}

\noindent {\bf Proof}.  Let $M=\sup\{\Re(h(z)), |z|<r\}$.  Note 
that $M\ge 0$.

\noindent Let $g(z) = {h(z)\over 2M-h(z)}$.  Then $|g(z)| \le 1$ 
if $|z|\le r$,
$g$ is analytic in $D_r$ and $g(0) = 0$.  By the Schwarz lemma 
we have
$$|g(z)|\le {|z|\over r}\ {\rm for\ all}\ z\ {\rm in}\ D_r\ {\rm 
and}\ |g'(0)|
\le 1/r.$$

\noindent Since $h(z) = {2M g(z)\over 1+g(z)}$ we have $h'(0) = 
2M g'(0)$ hence
$|h'(0)| \le 2M/r.\quad $\qed\vskip6pt

\noindent We now prove the main result of this section, 
namely that any weak
Hilbert space has the proportional UAP.  Let $X$ be a weak 
Hilbert
space.  We will show that $X$ satisfies (ii) in Proposition 1.2.  
Actually, we
obtain the following result of independent interest.\vskip6pt

\proclaim Theorem 1.5. {\sl Let $X$ be a weak Hilbert space 
with constant $C$ as
in (1.5).  Let $u,v$ be nuclear operators on $X$ and let $\rho$ 
be the spectral
radius of $u+v$.  Then if $\rk(u)\le n$ and $N(u)>1$
$$|\tr(u+v)|\le 2n \rho\ \Log N(u) + 2C + 2CN(v)\leqno (1.7)$$
hence also
$$N(u+v)\le 2n\|u+v\| \Log N(u) + 2C + 2CN(v).\leqno (1.8)$$}

\noindent {\bf Proof}:  Let $R=1/\rho$.  The function $f(z) = 
{\rm det}\,(I+z(u+v))$ is
entire and does not vanish in $D_R$.  Therefore there is an 
analytic function
$h$ on $D_R$ such that $f = \exp (h)$ and since $f(0)=1$ we 
can 
assume
$h(0)=0$.\vskip6pt

\noindent Note that $f'(0) = h'(0) = \tr(u+v)$.  By (1.6) we have 
if $N(u) \ge 1$
and $r<R$

$$\sup_{|z|=r} \Re\ h(z) \le n \Log N(u) + Cr + Cr N(v)$$

\noindent hence by Lemma 1.4

$$|\tr(u+v)| = |h'(0)|\le {2n\over r} \Log N(u) + 2C + 2CN(v)$$

\noindent Letting r tend to $R={1\over \rho}$, we obtain (1.7).

\noindent For (1.8) we simply observe that if $N(u)>1$ we have
$$N(u) = \sup \{|\tr(uS)|; \quad S:  X\rightarrow X,\quad 
\|S\|\le 
1,\quad
N(uS)>1\}.\leqno (1.9)$$

\noindent Therefore (1.8) follows from (1.7) since $\rho\le 
\|u+v\|$ and we can
take the supremum of (1.7) over all $S$ as in (1.9).  
\qed\vskip6pt

\noindent Finally we prove the ``only if'' part of our main 
theorem.  Let $X$ be
a weak Hilbert space.  The first and second authors proved,
respectively, that 
$X$ is reflexive (cf. [P 3], chapter 14) and that $X$ has the AP 
([P 3], chapter 15). We will show that (ii)
in Proposition 1.2 holds for suitable constants.  Let $T_1, T_2$ 
be as in
Proposition 1.2.  Let $u=T_1+T_2$ and $v=-T_1$.

\noindent By homogeneity we may assume 
$n\|T_2\| + N(T_1) = 1$.

\noindent Then if $N(T_1 + T_2) > 1$ we have by (1.8)
$$\eqalign{N(T_1+T_2)&\le N(u+v)+N(v)\cr
&\le 2n\|T_2\| \Log N(T_1+T_2) + 2C +  (2C+1)N(T_1)\cr 
&\le 2 \Log N(T_1 + T_2)+ 4C+1,}$$

\noindent and (since $2\Log x\le (x/2) + 2\ {\rm if}\ x>1$) this 
implies that if $N(T_1+T_2)>1$ then
$$ N(T_1+T_2)\le 8C+6.$$

\noindent Since in the case $N(T_1+T_2)\le 1$ this bound is 
trivial, we
conclude by homogeneity that (if $T_1+T_2$ has rank $\le n$)

$$N(T_1+T_2)\le (8C+6)(n\|T_2\| + N(T_1)).$$

\noindent By proposition 1.2 and 1.1 we conclude that $X$ has 
the proportional
UAP.\quad \qed\vskip6pt

\noindent {\bf Remark}:  Replacing $(u+v)$ by $\epsilon(u+v)$ 
in (1.7) yields 
that if
$\epsilon\ge N(u)^{-1}$, then $ |\tr(u+v)|\le 
2n\rho\Log(\epsilon
N(u))+2C\epsilon^{-1} + 2CN(v)$ hence after minimization over 
$\epsilon \ge
N(u)^{-1}$ we find that if $N(u)\ge n\rho/C$, then
$$|\tr(u+v)|\le 2n\rho(\Log^+({CN(u)\over n\rho})+1) 
+2CN(v).$$
\noindent On the other hand if $N(u)<n\rho/C$ we have 
trivially since $C\ge 1$
$$\eqalign{|\tr(u+v)|&\le N(u+v)\cr&\le n\rho/C+N(v)\cr&\le
2n\rho +2CN(v)\cr}$$
\noindent hence we conclude that without any restriction on 
$N(u)$ we have if
$\rk(u)\le n$
$$|\tr(u+v)|\le 2n\rho(\Log^+({CN(u)\over n\rho})+1)+2CN(v).$$

\noindent Even in the case of a Hilbert space we do not see a 
direct proof of this inequality.

\vfill\eject

\magnification\magstep1
\baselineskip=18pt
%lines for a solid square at end of proofs
\def \qed {\vrule height6pt  width6pt depth0pt}

\def \Bbb {\bf}
\def \dim {{\rm dim\,}}
\def \det {{\rm det\,}}
\def \Log{{\rm Log \,}}

\noindent {\bf $\S 2.$ The converse}\vskip6pt

\noindent Recall that $X$ is a weak cotype 2 space if there are 
constants $C$
and $0<\delta<1$ such that every finite dimensional subspace 
$E\subset X$
contains a subspace $F\subset E$ with $\dim F \ge \delta 
\,\dim E$ 
\ such that \ 
$d_F\equiv d(F, l^{\dim F}_2)\le C$  (cf. [MP]).\vskip6pt 
 \noindent We begin with a slightly modified presentation of  
Mascioni's [Ma 1] proof that a Banach  space $X$ which has
proportional UAP must have weak cotype $2$.  Suppose that 
$k_X(n, K) \le L\, n$
for all $n=1,2,\dots $.  Take any $(1+\delta)n$-dimensional 
subspace $G_0$ of
$X$, and, using Milman's subspace of quotient theorem [M] 
(or see [P 3], chapters 7 \& 8),  choose an $n$--dimensional 
subspace $G$ of $G_0$ for 
which there exists a
subspace $H$ of $G$ such that $\dim H \le \delta n$ and 
$d=d_{G/H}$ is bounded
by a constant which is independent of $n$, where $\delta$ is 
chosen so that
$\delta L \le {1\over 2}$. 

\noindent Take $T\colon X \to X$ with $T_{|H}=I_H$, $\Vert T 
\Vert \le K$,  and
$\rk(T) \le \delta Ln$. Let $Q\colon G \to G/H$ be the quotient 
map and  set 
$E= \ker(T)\cap G$.  If $x$ is in $E$, then 
$$\eqalign{\Vert Qx \Vert & = \, \inf_{h \in H} \, \Vert x-h 
\Vert   \cr &\ge \,
\inf_{h \in H} \, {\Vert (I-T) (x -h)\Vert \over \Vert I-T 
\Vert} \cr &\, =  {\Vert x  \Vert  \over \Vert I-T \Vert}\, ;}$$ 
that is, $Q_{|E}$ is an isomorphism
and \ $\Vert (Q_{|E})^{-1}\Vert \le \Vert I-T \Vert$. Thus \ 
$d_E \le \Vert I-T
\Vert d_{G/H}$, which finishes the proof since \ $\dim E \ge n 
- \delta Ln \ge
{n\over 2}$.\vskip6pt

\noindent Since we do not know {\it a priori} that the 
proportional UAP
dualizes, we need to prove Mascioni's theorem under a weaker 
hypothesis.

\proclaim Theorem 2.1.  {\sl Let $X$ be a Banach space.  
Assume that there are
constants $K$ and $L$ such that for all finite dimensional 
subspaces 
$E\subset X$ there is an
operator $T: X\rightarrow X$ satisfying $T_{|E} = I_E$ and such 
that $\|T\|\le
K$ and $\pi_2(T) \le L{(\dim E)}^{1/2}$.  Then $X$ is a weak 
cotype 2 space.}

\noindent {\bf Proof:} Recall that if \ $T: X\rightarrow Y$ \ is 
an operator,
then $$\pi_2(T)^2 = \sup \,\big{\{} \sum \Vert Tu(e_i) 
\Vert^2 ;  u: l_2
\rightarrow X, \,\, \Vert u \Vert \le 1\big{\}}$$ where \
$\{e_i\}_{i=1}^{\infty}$ \ is the unit vector basis for 
$l_2$.\vskip6pt
\noindent Take any $(1+\delta)n$-dimensional subspace $G_0$ 
of $X$, where $\delta$ is chosen so that $L\sqrt {\delta }\le 
{1\over 8}$.  By Milman's subspace of quotient theorem ([M] or 
[P 3], chapters 7 \& 8), we can choose an $n$--dimensional 
subspace $G$ of 
$G_0$ for which there exists a subspace $H$ of $G$ such that 
dim$(H) \le \delta n$ and $d=d_{G/H}$ is bounded by a 
constant which is independent of $n$.  Using the ellipsoid of 
maximal volume, we get from the Dvoretzky-Rogers lemma (cf. 
[P 3], lemma 4.13) a norm one operator
$J\colon l_2^n \to G$ such that $\Vert x_i \Vert \ge {1\over 
2}$ for all
$i=1,2,\dots ,{n\over 2}$, where $x_i=Je_i$.  By 
[FLM], all we need to
check is that there is a constant $\tau$ so that  
$$\left ( \hbox{Average}_{\pm}\, \Vert
\sum_{i=1}^n \pm x_i \Vert^2 \right )^{1\over 2} \; \ge \; \tau 
\sqrt n.$$

\noindent Let $Q\colon G \to G/H$ be the quotient map. Then
$$\tau_0  \sqrt {n/4} \; \equiv \left ( \hbox{Average}_{\pm}\, 
\Vert \sum_{i=1}^{n/2} \pm Qx_i \Vert^2 \right)^{1 \over 2} 
\ge \; {1 \over d} \left ( \sum_{i=1}^{n/2} \Vert Qx_i \Vert ^2 
\right )^{1 \over 2}.$$
So we can assume without loss of generality that $\Vert Qx_i 
\Vert \le d\tau_0$ for $1\le i \le {n\over 4}$.  Now take 
$T\colon X \to X$ with $T_{|H}=I_H$ and $\pi_2(T) \le L\sqrt 
{\delta n}$. Thus also
$$\left ( \sum_{i=1}^{n/4} \Vert Tx_i \Vert ^2 \right )^{1\over 
2} \le \; L\sqrt {\delta n}$$
hence without loss of generality $\; \Vert Tx_1 \Vert \, \le \, 
2L \sqrt {\delta }$.  But then
$$\eqalign{d\tau_0 \;\; &\ge \; \qquad \Vert Qx_1 \Vert \; 
\qquad = \;\; \inf_{h \in H} \, \Vert x_1-h \Vert   \cr \cr  & 
\ge \; \; {\Vert (I-T) x_1 \Vert \over \Vert I-T \Vert} \quad \; 
\ge \; \;  {{{1 \over 2} - \Vert Tx_1 \Vert } \over \Vert I-T 
\Vert} \cr \cr  & \ge \;\; {{{1\over 2} - {2L \sqrt {\delta}}} 
\over \Vert I-T \Vert} \;\;\;\; \ge \;\; {1\over {4 \Vert I-T 
\Vert}};}$$
that is, $\; \tau_0 \, \ge \, \left ( 4 d \Vert I-T \Vert \right )^{-
1}$. \hfill \qed\vskip6pt
\noindent {\bf Proof of converse of Main Theorem:}
By [LT 2], $X^{**}$ also has proportional UAP; in fact, \ 
$k_{X^{**}}(K,n) = k_X(K,n)$ for all \ $K$ \ and \ $n$.  Then, just 
as in the proof of Theorem 4 in [Ma 2], Lemma 1 in [Ma 2] (or 
its unpublished predecessor proved by Bourgain and 
mentioned in [Ma 2]) yields that  $X^{*}$   satisfies the 
hypothesis of Theorem 2.1 and hence  
$X^*$  as well as  $X$  has weak cotype 2.  By the results of 
[P 1], it only remains to check that   $X$  has non-trivial type; 
this is done as in Theorem 3.3 of [Ma 1]: since  $X$   and  $X^*$   
have 
weak cotype $2$, they both have cotype $2+\epsilon$ for all 
$\epsilon > 0$.  Since $X$ has the bounded approximation 
property, the main result of [P 1] yields that $X$  has 
non-trivial type. \qed\vskip6pt
\noindent{\bf Remark:} With a bit more work, the converse of 
the Main Theorem can be improved.  Following Mascioni [Ma 2], 
given an ideal norm $\alpha$, a normed space $X$, and $K>1$, 
we write
$$\alpha{\rm-}k_X(K,n) = \sup_{\scriptstyle E\subset 
X\atop\scriptstyle
\dim E=n }\ \ \inf\ \{\alpha(T)\}$$
\noindent where the infimum runs over all finite rank 
operators $T: X\rightarrow X$ such that $\|T\|
< K$ and $T_{|E} = I_{|E}$. (We use ``$<K$" instead of  ``$\le K$" 
in order to avoid in the sequel statements involving awkward 
``$K+\epsilon$ for all $\epsilon >0$".) We say that $X$ has the 
$\alpha$-UAP if there is a $K>1$ such that for all $n$,  
$\alpha$-$k_X(K,n) < \infty$; when the value of $K$ is 
important, we say 
that   $X$ has the $K$-$\alpha$-UAP.  Notice that the ``finite 
rank" can be ignored if the space $X$ has the metric 
approximation property or (by adjusting $K$) if $X$ has the 
bounded approximation property.  Here we are interested in 
$\alpha=\pi_2$ and $\alpha=\pi_2^d$, where $\pi_2^d(T) 
\equiv\pi_2(T^*)$.  Since for either of these $\alpha$'s, 
$\alpha(T) < \infty$ implies that $T^2$ is uniformly 
approximable by finite rank operators ($T^2$ factors through a 
Hilbert-Schmidt operator), for these two $\alpha$'s the 
K-$\alpha$-UAP implies the bounded approximation property.  
In fact, by passing to ultraproducts and using [H], it follows that 
for either of these $\alpha$'s the K-$\alpha$-UAP implies the 
($K^2+\epsilon$)-UAP; in particular, $X^{**}$ has the bounded 
approximation property.  (This is really a sloppy version of 
Mascioni's reasoning 
[Ma 2]; Mascioni gives a better estimate for $k_X(K',n)$ in 
terms of $\alpha$-$k_X(K,n)$.)  We now state an improvement 
of 
the converse in the Main Theorem: 
\vfill \eject
\proclaim Theorem 2.2. {\sl Suppose that there are constants 
$K$ and $L$ so that the Banach space $X$ satisfies for all $n$
\ $\pi_2$-$k_X(K,n) \le L\sqrt n$ \ and \ 
$\pi_2^d$-$k_X(K,n) \le L\sqrt n$.  Then $X$ is a weak Hilbert 
space.}

\noindent {\bf Proof:} In view of the discussion above, we can 
ignore the ``finite rank" condition in the definition of 
$\alpha$-UAP.  It is then easy to see for $\alpha=\pi_2$ or 
$\alpha=\pi_2^d$ that for all $n$ and $K$
$\alpha$-$k_{X^{**}}(K,n) = \alpha$-$k_X(K,n)$, hence by 
Lemma 1 
of [Ma 2] and Theorem 2.1 we conclude that $X$ and $X^*$ have 
weak cotype 2.  The argument used in the proof of the 
converse in the Main Theorem shows that $X$ has non-trivial 
type, so $X$ is weak Hilbert. \qed

\vfil\eject

\vfill\eject

\noindent {\bf $\S 3.$ Related results and concluding 
remarks}\vskip6pt
\noindent In [Ma 1], Mascioni proved (but stated in slightly weaker form)
that
for $2<p<\infty$ and  all $K$, there exists $\delta=\delta(p,K)>0$ so that for
all $n$,  $k_p(K,n) \ge \delta  n^{p/2}$. (We write $k_p(K,n)$ for 
$k_{L_p}(K,n)$ and 
$\alpha$-$k_p(K,n)$ for $\alpha$-$k_{L_p}(K,n)$.  See [FJS], 
[JS], and [Ma 1] for results about the UAP in $L_p$-spaces, 
$1\le p \le \infty$.) We prove 
here the corresponding result for $1 < p < 2$. 
\proclaim Theorem 3.1. {\sl For each $2<p<\infty$ and $K>1$, 
there exists $\epsilon=\epsilon(p,K)>0$ so that for all $n$, 
$\pi_2$-$k_p(K,n)  \ge \epsilon  n^{p\over 4}$. 
Consequently, for $1<q<2$,  
$k_q(K,n) \ge \pi_2^d$-$k_q(K,n)^2  \ge \epsilon^2 n^{p\over 
2}$  
where  ${1\over p} + {1\over q} = 1$. }

\noindent {\bf Proof:} The proof is basically the same as the 
proof of Theorem 2.1 once we substitute a result of Gluskin [Gl] 
for Milman's
subspace of quotient theorem, so we use notation similar to 
that used in Theorem
2.1.  Fix $n$, set $G=l_p^n$, let $J$ denote the formal identity 
from $l_2^n$
into $l_p^n$, and let $x_i = Je_i$ be the unit vector basis of 
$l_p^n$. By
Gluskin's theorem [Gl], given any $\gamma>0$ there is 
$M=M(p,\gamma)$
independent of $n$ and a subspace $H$ of $G$ with $\dim H \le 
Mn^{2\over p}$ so
that  $d\equiv d_{G/H} \le \gamma n^{{1\over 2}-{1\over p}}$.  
Let $Q\colon G \to G/H$ be the quotient map.  
Define $\epsilon_0$ by the formula   
$\pi_2$-$k_p(K,\dim H) = \epsilon_0  (\dim H)^{p\over 4}$ 
and choose  $T\colon X \to X$ with $T_{|H}=I_H$,
$\Vert T \Vert \le K$,  and 
$\pi_2(T) \le \epsilon_0  M^{p\over 4}\sqrt{n}$.  
We need to show that
$\epsilon_0$ is bounded  away from
\ $0$ \ independently of $n$.  Now
$$n^{1\over p} = \left ( \hbox{Average}_{\pm}\, 
\Vert \sum_{i=1}^n \pm x_i \Vert^2 \right)^{1 \over 2} 
\ge \;\left ( \hbox{Average}_{\pm}\, 
\Vert \sum_{i=1}^n \pm Qx_i \Vert^2 \right)^{1 \over 2} 
\ge \; {1 \over d} \left ( \sum_{i=1}^n \Vert Qx_i \Vert ^2 
\right )^{1 \over 2},$$
So we can assume without loss of generality that for 
$i=1,\dots {n\over2}$,
$$\Vert Qx_i \Vert \le \sqrt 2 d n^{{1\over p}-{1\over 2}}.$$

\noindent Since \ 
$\left ( \sum_{i=1}^{n/2} \Vert Tx_i \Vert ^2 \right )^
{1\over 
2} \le \; \pi_2(T) \le \epsilon_0  M^{p\over 4}\sqrt{n}$, 
we can also assume without loss of generality that 
$\; \Vert Tx_1 \Vert \, \le \, \sqrt{2} \epsilon_0 M^{p\over 
4}$. 
But then 
$$\eqalign{\sqrt 2 d n^{{1\over p}-{1\over 2}}\;\; &\ge \; 
\qquad \Vert Qx_1 \Vert \; 
\qquad = \;\; \inf_{h \in H} \, \Vert x_1-h \Vert   \cr \cr  & 
\ge \; \; {{\Vert (I-T) x_1 \Vert} \over {\Vert I-T \Vert}} 
\quad \; 
\ge \; \;  {{1 - \Vert Tx_1 \Vert } \over {\Vert I-T 
\Vert}} \cr \cr  & \ge \;\; {{1 - \sqrt{2} \epsilon_0 M^{p\over 
4}} 
\over {\Vert I-T \Vert}} \; ;}$$
that is,
$$d \; \ge \; {{1-\sqrt{2} \epsilon_0 M^{p\over 
4}}\over{\sqrt{2}\Vert I-T \Vert}}n^{{1\over 2}-{1\over p}}.$$
For sufficiently small \ $\gamma$ (e.g.,  
$\gamma \le {1\over 2}(K+1)^{-1}$), this 
 gives a lower bound on $\epsilon_0$ since 
$d \le \gamma n^{{1\over 2}-{1\over p}}$.\vskip6pt

\noindent The ``consequently" statement follows by duality 
from Lemma 1 in [Ma 2]. \qed \vskip6pt

{\centerline {*\quad\quad*\quad\quad*}}\vskip6pt

\noindent The trick of Mascioni's [Ma 1] mentioned at the 
beginning of section 2 can be used to answer a question 
Pe\l czy\'nski asked the authors twelve years ago; namely, 
whether 
every $n$-dimensional subspace of $l_\infty^{2n}$ 
well-embeds into $l_\infty^{(1+\epsilon)n}$ for each 
$\epsilon > 0$.  
Since $l_\infty^{2n}$ has an $n$-dimensional quotient $F$ with 
$d_F$ bounded independently of $n$ by Ka\v sin's theorem
([K] or [P 3], corollary 6.4), Proposition 
3.2 gives a strong negative answer to Pe\l czy\'nski's 
question.\vskip6pt

\proclaim Proposition 3.2. {\sl Set $G=l_\infty^n$, let $H$ be a 
subspace of $G$, set $d=d_{G/H}$, and assume that $H$ is 
$K$-isomorphic to a subspace of $l_\infty^k$.  Then  $d \ge 
e^{-2}(K+1)^{-1}\bigl({{n-k}\over {\Log k}}\bigr)^{1\over 2}$.}

\noindent {\bf Proof:} Let $u\colon H \to l_\infty^k$ satisfy 
$\Vert u \Vert =1$, $\Vert u^{-1}_{|uH} \Vert \le K $, let $U$ 
be a 
norm one extension of $u$ to an operator from $G$ to 
$l_\infty^k$, let $S$ be an extension of $u^{-1}$ to an operator 
from $l_\infty^k$ to $G$ with 
$\Vert S \Vert \le K $, and set $T = SU$.  So $T_{|H} = I_H$ and 
$\Vert T \Vert \le K $.  Let $Q \colon G \to G/H$ be the 
quotient map and set  $E = \ker(T)$, so that 
$\dim E \ge n-k$.  Thus (see the argument at the beginning of 
section 2)
$$ d_E \le \Vert I-T \Vert d \le (K+1)d\, .$$
But by [BDGJN], p. 182 (let $s=\Log k$ there), 
$$d_E \ge e^{-2}\Bigl({{n-k}\over {\Log k}}\Bigr)^{1\over 2}.\ \ 
\qed $$
\vskip6pt

{\centerline {*\quad\quad*\quad\quad*}}\vskip6pt

\noindent We conclude with two open problems related to the 
material in section 1.  \vskip6pt
\noindent {\bf Problem 3.3.} {\sl If 
$X$ is a weak Hilbert space, then is 
$k_X(K,n)$ proportional to $n$ for all $K>1$?}\vskip6pt
\noindent Since weak Hilbert spaces are superreflexive, 
for all $K>1$ 
$k_X(K,n) < \infty$ for every weak Hilbert space $X$ by a 
result of Lindenstrauss and Tzafriri [LT 2], but their argument 
does not give a good estimate of $k_X(K,n)$ for $K$ close to one 
when one has a good estimate for large $K$.\vskip6pt

\noindent For the known weak Hilbert spaces $X$, the growth 
rate of \ 
$k_X(K,n)-n$ \ is very slow (cf. [J], [CJT]), at least for 
sufficiently 
large $K$.  It follows from recent work of Nielsen and 
Tomczak-Jaegermann that \ $k_X(K,n)-n$ \ has the same kind 
of slow 
growth for any weak Hilbert space $X$ which has an 
unconditional basis. On the other hand, we do not know any 
improvement of the result presented in section 2 for general 
weak Hilbert spaces.  This suggests:\vskip6pt
\noindent {\bf Problem 3.4.} {\sl If $X$ is a weak Hilbert space, 
does there exist $K$ so that $k_X(K,n)-n = o(n)$?}

\vfill
\eject

\centerline {\bf References}\vskip6pt

\item{[BDGJN]}  G. Bennett, L.~E. Dor, V. Goodman, W.~B. 
Johnson, and C.~M. 
Newman, {\sl  On uncomplemented subspaces of  $L_p$, $1 < p 
< 2$,} {\bf 
Israel J. Math. 26} (1977), 178--187.\vskip6pt

\item{[CJT]} P.~G. Casazza, W.~B. Johnson, and L. Tzafriri,
{\sl On Tsirelson's space,} {\bf Israel J. Math. 47} (1984), 81--
98.\vskip6pt

\item{[FJS]}  T. Figiel, W. B. Johnson, and G. Schechtman, {\sl   
Factorizations of natural embeddings of  $l_p^n$ into  $L_r$, I,} 
{\bf   
Studia Math. 89} (1988), 79--103.\vskip6pt

\item {[FLM]} T. Figiel, J. Lindenstrauss, and V. Milman,  {\sl 
The dimension of
almost spherical sections of convex bodies,}  {\bf Acta Math. 
139} (1977), 53--94.\vskip6pt

\item{[Gl]} E. D. Gluskin, {\sl Norms of random matrices and 
diameters of finite dimensional sets,} {\bf Mat. Sbornik 120} 
(1983), 180--189.\vskip6pt

\item {[G]}  A. Grothendieck,  {\sl La Th\'eorie de Fredholm,}  
{\bf Bull. Soc. Math.
France 84} (1956), 319--384.\vskip6pt

\item{[H]} S. Heinrich, {\sl Ultraproducts in Banach space 
theory,} {\bf J. Reine Angew. Math. 313} (1980), 72--
104.\vskip6pt

\item{[J]} W. B. Johnson, {\sl Banach spaces all of whose 
subspaces have  
the approximation property,} {\bf  Special Topics of Applied  
Mathematics} 
North-Holland (1980), 15--26.\vskip6pt

\item{[JS]}  W. B. Johnson and G. Schechtman, {\sl Sums of 
independent  
random variables in rearrangement invariant function spaces,} 
{\bf  Ann.  Prob. 17} (1989), 789--808.\vskip6pt

\item{[K]} B. S. Ka\v sin, {\sl Sections of some 
finite-dimensional sets and  
classes of smooth functions,} {\bf  Izv. Akad. Nauk SSSR 41} 
(1977),  334--351 (Russian).\vskip6pt

\item {[LT 1]}  J. Lindenstrauss and L. Tzafriri,  {\sl On the 
complemented subspaces
problem,}  {\bf Israel J. Math. 9} (1971), 263--269.\vskip6pt

\item {[LT 2]}  J. Lindenstrauss and L. Tzafriri,  {\sl The 
uniform approximation property in Orlicz spaces,} {\bf Israel J. 
Math. 23} (1976), 142--155.\vskip6pt

\item{[Ma 1]} V. Mascioni, {\sl Some remarks on the uniform 
approximation 
property in Banach spaces,} {\bf Studia Math.} (to 
appear).\vskip6pt

\item {[Ma 2]}  V. Mascioni, {\sl On the duality of the uniform 
approximation
property in Banach spaces,}  {\bf Illinois J. Math.} (to 
appear)\vskip6pt

\item{[M]} V. D. Milman, {\sl Almost Euclidean quotient spaces 
of  
subspaces of finite dimensional normed spaces,} {\bf  Proc. 
AMS 94}  
(1985), 445--449. \vskip6pt

\item {[MP]}  V. Milman and G. Pisier,  {\sl Banach spaces with 
a 
weak cotype 2
property,}  {\bf Israel J. Math. 54} (1986), 139--158.\vskip6pt

\item {[PT]}  A. Pajor and N. Tomczak-Jaegermann,  {\sl 
Subspaces of small
codimension of finite dimensional Banach spaces,}  {\bf Proc. 
A.M.S. 97} (1986), 637--642.\vskip6pt

\item{[P 1]} G. Pisier, {\sl On the duality between type and 
cotype,} {\bf Springer Lecture Notes 939} (1982), 131--
144.\vskip6pt

\item {[P 2]}  G. Pisier,  {\sl Weak Hilbert spaces,}  {\bf Proc. 
London Math. Soc. 56}
(1988), 547--579.\vskip6pt

\item {[P 3]}  G. Pisier,  {\sl The volume of convex bodies and 
Banach space
Geometry,}  {\bf Cambridge Univ. Press} (1989).\vskip6pt

\item {[S]}  A. Szankowski, {\sl On the uniform approximation
property in Banach spaces,}  {\bf Israel J. Math. 49} (1984), 
343--359. \vskip20pt

\noindent Texas A\&M University, College Station, TX 77843, 
U.S.A.

\noindent Texas A\&M University, College Station, TX 77843, 
U.S.A., and Equipe d'Analyse, Universit\'e Paris VI, 75230 
Paris, FRANCE

\vfill\eject

\end